\documentclass[12pt]{amsart}

\usepackage{amssymb}
\usepackage{amsmath}
\usepackage{graphicx}

\newtheorem{theorem}{Theorem}
\newtheorem{lemma}[theorem]{Lemma}
\newtheorem{proposition}[theorem]{Proposition}
\newtheorem{corollary}[theorem]{Corollary}
\newtheorem{_definition}[theorem]{Definition}

\newtheorem{_remark}[theorem]{\it Remark}
\newenvironment{remark}{\begin{_remark}\rm}{\end{_remark}}

\newcommand{\F}{\mathord{\mathbb F}}
\renewcommand{\P}{\mathord{\mathbb  P}}
\newcommand{\Z}{\mathord{\mathbb Z}}

\newcommand{\HHH}{\mathord{\mathcal H}}
\newcommand{\III}{\mathord{\mathcal I}}
\newcommand{\PPP}{\mathord{\mathcal P}}
\newcommand{\RRR}{\mathord{\mathcal R}}

\newcommand{\VVV}{\mathord{\mathcal V}}

\newcommand{\GL}{\mathord{\mathrm {GL}}}
\newcommand{\PGU}{\mathord{\mathrm {PGU}}}
\newcommand{\PGL}{\mathord{\mathrm{PGL}}}

\newcommand{\Aut}{\operatorname{\mathrm {Aut}}\nolimits}
\newcommand{\Stab}{\mathord{\mathrm{Stab}}}

\newcommand{\inv}{\sp{-1}}

\newcommand{\sptimes}{\sp{\times}}

\newcommand{\Fqsq}{\F_{q^2}}
\newcommand{\tlG}{\widetilde{G}}
\newcommand{\Pn}{\P^n}
\newcommand{\acF}{k}

\newcommand{\transpose}[1]{\hskip 1pt {}^t\hskip -.4pt{#1}}

\newcommand{\inj}{\hookrightarrow}
\newcommand{\isom}{\mathbin{\,\raise -.6pt\rlap{$\to$}\raise 3.5pt \hbox{\hskip .3pt$\mathord{\sim}$}\,}}
\newcommand{\maprightsp}[1]{\; \smash{\mathop{\; \longrightarrow \; }\limits\sp{#1}}\; }
\newcommand{\mapdown}{\phantom{\Big\downarrow}\hskip -8pt \downarrow}

\newcommand{\set}[2]{\{\; {#1} \; \mid \; {#2} \;  \}}

\newcommand{\rmand}{\textrm{and}}
\newcommand{\quand}{\quad\rmand\quad}

\numberwithin{equation}{section}

\begin{document}

\title[Rational normal curves and a Hermitian variety]%
{A  note on rational normal curves totally tangent to a Hermitian variety
}

\author{Ichiro Shimada}
\address{
Department of Mathematics, 
Graduate School of Science, 
Hiroshima University,
1-3-1 Kagamiyama, 
Higashi-Hiroshima, 
739-8526 JAPAN
}
\email{shimada@math.sci.hiroshima-u.ac.jp
}

\thanks{Partially supported by
 JSPS Grants-in-Aid for Scientific Research (B) No.20340002 
}



\begin{abstract}
Let $q$ be a power of a  prime integer $p$,
and let $X$ be a  Hermitian variety of degree $q+1$ in the $n$-dimensional
projective space.
We count the number of rational normal curves that are tangent to $X$
at distinct $q+1$ points with intersection multiplicity $n$.
This generalizes a result of B.~Segre on 
the permutable pairs of a Hermitian curve and a smooth conic.
\end{abstract}

\maketitle

%
%

\section{Introduction}
Throughout this paper,
we fix a power $q:=p^\nu$ of a  prime integer $p$.
Let $\acF$ denote the algebraic closure of the finite field $\Fqsq$.
\par
\medskip
Let $n$ be an integer $\ge 2$.
We say that a hypersurface $X$ of $\Pn$
defined over $\Fqsq$ 
is a \emph{Hermitian variety}
if $X$ is  
projectively isomorphic  over $\Fqsq$ 
to the Fermat variety
$$
X_I:=\{x_0^{q+1}+\cdots+ x_{n}^{q+1}=0\}\;\subset\; \Pn
$$
of degree $q+1$.
(Strictly speaking, one should say that $X$ is
a Hermitian variety of rank $n+1$. 
Since we treat only nonsingular Hermitian varieties in this paper,
we omit  the term ``of rank $n+1$".)
We say that a hypersurface $X$ of $\Pn$ 
defined over $\acF$
is a \emph{$\acF$-Hermitian variety}
if $X$ is  projectively isomorphic  over $\acF$ 
to $X_I$.
By definition, the projective automorphism group
$\Aut(X)\subset \PGL_{n+1}(\acF)$ of a  $\acF$-Hermitian variety $X$ is 
conjugate  to $\Aut(X_I)=\PGU_{n+1}(\Fqsq)$ in $\PGL_{n+1}(\acF)$.
\par
\medskip
Let $X$ be a  $\acF$-Hermitian variety in $\P^n$.
A rational normal curve $\Gamma$ in $\Pn$ defined over $\acF$
is said to be \emph{totally tangent to $X$}
if $\Gamma$ is tangent to $X$ at distinct $q+1$ points
and the intersection multiplicity at each intersection point is $n$.
\par
\medskip
A subset $S$ of a rational normal curve $\Gamma$ is called a \emph{Baer subset}
if there exists a coordinate $t: \Gamma\isom \P^1$ on $\Gamma$ such that
$S$ is the inverse image  by $t$ of the set
$\P^1(\F_q)=\F_q\cup\{\infty\}$ of  $\F_q$-rational points of $\P^1$.
\par
\medskip
The purpose of this paper is to prove the following:
\begin{theorem}\label{thm:main}
Suppose that $n\not\equiv 0 \pmod p$ and $2n\le q$.
Let $X$ be a  $\acF$-Hermitian variety in $\P^n$.

{\rm (1)}
The set $R_X$ of rational normal curves totally tangent to $X$
is non-empty, and 
$\Aut(X)$ acts on $R_X$ transitively
with the stabilizer subgroup isomorphic to $\PGL_{2}(\F_q)$.
In particular, we have
$$
|R_X|=|\PGU_{n+1}(\Fqsq)|/|\PGL_{2}(\F_q)|.
$$

{\rm (2)}
For any  $\Gamma\in R_X$, the points in 
$\Gamma\cap X$ form a Baer subset of $\Gamma$.

{\rm (3)}
If $X$ is a Hermitian variety,
then every $\Gamma\in R_X$ is defined over $\Fqsq$ 
and every point of $\Gamma\cap X$ is $\Fqsq$-rational.
\end{theorem}
The study of Hermitian varieties was initiated by B.~Segre in~\cite{MR0213949}.
Since then, Hermitian varieties have been intensively studied mainly from
combinatorial point of view in finite geometry. (See, for example, Chapter 23 of~\cite{MR1363259}). 
B.~Segre obtained  Theorem~\ref{thm:main} for the case $n=2$
in the investigation of commutative pairs of polarities~\cite[n.~81]{MR0213949}. 
We give a simple proof of the higher-dimensional analogue
(Theorem~\ref{thm:main}) of his result using arguments of projective geometry over $\acF$.
\par
\medskip
{\bf Notation.}
(1) For simplicity, we put
$$
\tlG:=\GL_{n+1}(\acF)
\quand
G:=\PGL_{n+1}(\acF).
$$
We let $G$ act on $\Pn$ from \emph{right}.
For $T\in \tlG$,
we denote by $[T]\in G$ the image of $T$ 
by the natural homomorphism $\tlG\to G$.
The entries $a_{i,j}$ of a matrix $A=(a_{i,j})\in \tlG$ are indexed by
$$
N:=\set{(i, j)\in \Z^2}{0\le i\le n,\;\; 0\le j\le n}.
$$
\par
(2)
Let  $M$ be a matrix  with entries in $\acF$.
We denote by $\transpose{M}$ the transpose of $M$,  and 
 by $\overline{M}$ the matrix obtained from $M$ by applying 
 $a\mapsto a^q$ to the entries.
\section{Proof of Theorem~\ref{thm:main}}
The following well-known proposition
is the main tool of the proof.
\begin{proposition}
The map
$\lambda : \tlG \to \tlG$
defined by $\lambda(T):=T\transpose{\overline{T}}$
is surjective.
The image of $\GL_{n+1}(\Fqsq)\subset  \tlG$
by $\lambda$ is equal to the set  
$$
\HHH:=\set{H\in \GL_{n+1}(\Fqsq)}{H=\transpose{\overline{H}}}
$$
of Hermitian matrices over $\Fqsq$.
\end{proposition}
\begin{proof}
The first part is a variant of Lang's theorem~(see~\cite{MR0086367} or~\cite{MR0466336}),
and can be proved by means of differentials~(see~\cite[16.4]{MR1102012}).
The second part is due to B.~Segre~\cite[n.~3]{MR0213949}.
See also~\cite[Section 1]{MR1794260}.
\end{proof}
For a matrix $A=(a_{i,j})\in \tlG$, we define  a homogeneous polynomial $f_A$ of  degree $q+1$ by 
$$
f_A:=\sum_{(i, j)\in N} a_{i,j}x_i x_j^q.
$$
If $g=[T]\in G$, 
then the image $X_I^{g}$ of the Fermat hypersurface $X_I$ by $g$ is defined by $f_{\lambda(T\inv)}=0$.
Hence we obtain the following:
\begin{corollary}
A  hypersurface $X$ of $\Pn$ is a  $\acF$-Hermitian variety
if and only if there exists a matrix $A\in \tlG$ such that
$X$ is defined by $f_A=0$. 
A  hypersurface $X$  is a  Hermitian variety
if and only if there exists 
a Hermitian matrix $H\in \HHH$ over $\Fqsq$ such that $X$ is defined by $f_H=0$.
\end{corollary}
Let $\VVV$ denote the set of all  $\acF$-Hermitian varieties in $\P^n$.
For $A\in \tlG$,
let $X_A\in \VVV$ denote the hypersurface defined by $f_A=0$.
(The Fermat hypersurface $X_I$ is defined by $f_I=0$,
where $I\in \tlG$ is the identity matrix.)
Remark that $G$ acts on $\VVV$ transitively.
\par
\medskip
Let $\RRR$ denote the set of all rational normal curves in $\P^n$,
and let $\PPP$ be the set of  pairs
$$
[\Gamma, (Q_0, Q_1, Q_\infty)],
$$
where $\Gamma\in \RRR$, and $Q_0, Q_1, Q_\infty$ are ordered three distinct points of $\Gamma$.
Let $\Gamma_0\in \RRR$ be the image of the morphism
$\phi_0: \P^1\inj \P^n$ given by
$$
\phi_0 (t):=[1:t:\dots: t^n]\in \P^n.
$$
We put 
$$
P_0:=\phi_0(0), \quad P_1:=\phi_0(1), \quad P_\infty:=\phi_0(\infty).
$$ 
Then we have $[\Gamma_0, (P_0, P_1, P_\infty)]\in \PPP$.
\begin{lemma}
The action of $G$ on $\PPP$ is simply transitive.
\end{lemma}
\begin{proof}
The action of $G$ on $\RRR$ is transitive by the definition of 
rational normal curves.
Let $\Sigma_0\subset G$ denote the stabilizer subgroup of $\Gamma_0\in \RRR$.
Then we have a natural homomorphism 
$$
\psi: \Sigma_0 \to \Aut(\Gamma_0)\cong \PGL_{2}(\acF).
$$
Note that $\Aut(\Gamma_0)$ acts on the set of ordered three distinct points of $\Gamma_0$ simple-transitively.
Hence
it is enough to show that $\psi$ is an isomorphism.
Since $\Gamma_0$ contains $n+2$ points such that any $n+1$ of them are linearly independent in $\P^n$,
$\psi$ is injective.
Since $\PGL_{2}(\acF)$ is generated by the linear transformations 
$$
t\mapsto at +b \quand t\mapsto 1/(t-c), \quad\textrm{where}\quad
a\in \acF\sptimes \;\;\textrm{and}\;\; b, c\in \acF, 
$$
it is enough to find matrices $M_{a, b}\in \tlG$ and $N_c\in \tlG$ such that
\begin{eqnarray*}
{}[1, at +b, \dots, (at +b)^n]&=& [1, t, \dots, t^n]\, M_{a, b} \quand\\ 
{}[(t-c)^n, (t-c)^{n-1}, \dots, 1]&=& [1, t, \dots, t^n]\, N_c 
\end{eqnarray*}
hold for any $a\in \acF\sptimes$ and $b, c\in \acF$.  This is immediate.
\end{proof}
%
We denote by $\III\subset \VVV\times \PPP$ the set  of all triples
$[X, \Gamma, (Q_0, Q_1, Q_\infty)]$ such that
\begin{enumerate}
\item $X\in \VVV$ and $[\Gamma, (Q_0, Q_1, Q_\infty)]\in \PPP$,
\item $\Gamma$ is totally tangent to $X$, and
\item $Q_0, Q_1, Q_\infty$ are contained in $\Gamma\cap X$.
\end{enumerate}
We then  consider the incidence diagram
$$
\begin{array}{ccc}
\III &\maprightsp{p_1} & \VVV \\
 {\raise 2pt \hbox to 0pt {\hskip -12pt \tiny $p_{23}$}}\mapdown && \\
\PPP. &&
\end{array}
$$
Note that $G$ acts on $\III$, and that the projections $p_1$ and $p_{23}$ are $G$-equivariant.
\par
\medskip
We consider the Hermitian matrix $B=(b_{i,j})\in \HHH$, where
$$
b_{i,j}:=\begin{cases}
\binom{n}{i}(-1)^i& \textrm{if $i+j=n$, }\\
0 & \textrm{otherwise.}
\end{cases}
$$
Then we have
$$
\phi_0^* f_B=\sum_{i=0}^n \binom{n}{i}(-1)^i t^i t^{q(n-i)} =(t^q-t)^n, 
$$
and hence 
$[X_B, \Gamma_0, (P_0, P_1, P_\infty)]$ is a point of $\III$.
Since $\III\ne \emptyset$ and the action of $G$ on $\VVV$ is transitive,
the $G$-equivariant map $p_1$ is surjective.
Thus $R_X\ne \emptyset$ holds for any $X\in \VVV$.
\par
\medskip
The following proposition is proved in the next section.
\begin{proposition}\label{prop:main}
Suppose that $n\not\equiv 0\pmod p$ and $n\le 2q$.
Then the  fiber of  $p_{23}$ over $[\Gamma_0, (P_0, P_1, P_\infty)]\in \PPP$
consists of a single point $[X_B, \Gamma_0, (P_0, P_1, P_\infty)]\in \III$.
In particular, $p_{23}$ is a bijection.
\end{proposition}
Theorem~\ref{thm:main} follows from Proposition~\ref{prop:main}
as follows.
First note that, for any  $X\in \VVV$, the map
$[X, \Gamma, (Q_0, Q_1, Q_\infty)]\mapsto \Gamma$
gives a surjection 
\begin{equation*}\label{eq:surjfib}
\rho_X \;:\; p_1\inv (X)\;\to\; R_X.
\end{equation*}
Proposition~\ref{prop:main} implies that $G$ acts on $\III$ simple-transitively.
If $S\subset \Gamma$ is a Baer subset of $\Gamma\in \RRR$,
then $S^g\subset \Gamma^g$ is a Baer subset of $\Gamma^g$ for any $g\in G$.
Since $\Gamma_0\cap X_B$ is a Baer subset of $\Gamma_0$,
we see that $\Gamma\cap X$ is a Baer subset of $\Gamma$
for any $[X, \Gamma, (Q_0, Q_1, Q_\infty)]\in \III$.
Therefore the assertion (2) follows.
Since $p_1$ is $G$-equivariant,
the stabilizer subgroup $\Aut(X)$ of $X$ in $G$ acts on the fiber $p_1\inv (X)$
simple-transitively
for any $X\in \VVV$.
Note that $\rho_X$ is $\Aut(X)$-equivariant.
Hence 
the stabilizer subgroup $\Stab(\Gamma)$ of $\Gamma\in R_X$ in $\Aut(X)$
acts on the fiber $\rho_X\inv(\Gamma)$ 
simple-transitively.
Moreover, since $\Gamma$ contains $n+2$ points such that any $n+1$ of them are linearly independent,
$\Stab(\Gamma)$ is embedded into $\Aut(\Gamma)\cong \PGL_{2}(\acF)$.
Since $\rho_X\inv(\Gamma)$ is 
the set of ordered three distinct points of the Baer subset
$\Gamma\cap X$ of $\Gamma$,
we see that
$\Stab(\Gamma)$ is conjugate to $\PGL_{2}(\F_q)$
as a subgroup of $\Aut(\Gamma)\cong \PGL_{2}(\acF)$.
Thus the assertion (1) follows.
The assertion (3) is immediate from the facts that $X_B$ is Hermitian,
that $\Gamma_0$ is defined over $\Fqsq$,
and that every points of $\Gamma_0\cap X_B$ is  $\Fqsq$-rational.
\qed
\section{Proof of Proposition~\ref{prop:main}}
Suppose that $A=(a_{i,j})\in \tlG$ satisfies 
$[X_A, \Gamma_0, (P_0,P_1, P_\infty)]\in \III$.
We will show that $A=c\, B$
for some $c\in \acF\sptimes$.
\par
\medskip
By the definition of $\III$,
there exists a polynomial $h\in \acF[t]$
such  
that the polynomial
$$
\phi_0^*f_A=\sum_{(i, j)\in N} a_{i,j} t^{i+qj}
$$
is equal to $h^n$,
 and that,
 regarded as a polynomial of degree $q+1$,
 $h$ has distinct $q+1$ roots including $0$, $1$ and $\infty$.
 In particular, we have $\deg h=q$ and $h(0)=0$.
Thus we can set
 $$
 h=\sum_{\nu=1}^q b_{\nu} t^{\nu}.
 $$
 Since $\deg h=q$ and since $t=0$ is a simple root of $h=0$,
we have 
\begin{equation*}\label{eq:bqnonzero}
b_q\ne 0\quand b_1\ne 0.
\end{equation*}
Let $c_m$ denote the coefficient of $t^m$ in $\phi_0^*f_A$.
We have $c_m=0$ if no $(i, j)\in N$  satisfy $i+qj=m$.
By the assumption $2n\le q$,
we have
\begin{equation}\label{eq:cm}
c_m=0\quad\textrm{if}\quad n<m<q\;\;\textrm{or}\;\; n+q(n-1)<m<qn.
\end{equation}
%
%
%
We will show that
\begin{equation}\label{eq:bzero1}
b_{\mu}= 0\quad\textrm{if}\quad n<\mu<q.
\end{equation}
Let $l$ be the largest integer such that $l<q$ and $b_l\ne 0$.
Since $n\not\equiv 0\pmod p$ and $b_q\ne 0$,
the coefficient $n\,b_q^{n-1} b_l$ of $t^{l+q(n-1)}$ in $h^n$ is non-zero.
Therefore $c_{l+q(n-1)}\ne 0$ follows from $\phi_0^*f_A=h^n$.
By~\eqref{eq:cm} and $l<q$, we have $l\le n$.
Hence~\eqref{eq:bzero1} holds.
In the same way, 
we will show that
\begin{equation}\label{eq:bzero2}
b_{\mu}= 0\quad\textrm{if}\quad 1<\mu < q-n+1.
\end{equation}
Let $l$ be the smallest integer such that $l>1$ and $b_l\ne 0$.
Since $n\not\equiv 0\pmod p$ and $b_1\ne 0$,
the coefficient $n\,b_1^{n-1} b_l$ of $t^{n-1+l}$ in $h^n$ is non-zero.
Therefore $c_{n-1+l}\ne 0$ follows from $\phi_0^*f_A=h^n$.
By~\eqref{eq:cm} and $l>1$, we have $n-1+l\ge q$.
Hence~\eqref{eq:bzero2} holds.
\par
\medskip
Combining~\eqref{eq:bzero1},~\eqref{eq:bzero2} with the assumption $2n\le q$,
we see that $h$ is of the form $b_q t^q+b_1 t$.
Since $h(1)=0$, we have
$$
h=b(t^q-t)\quad\textrm{for some $b\in \acF\sptimes$}.
$$
From $\phi_0^*f_A=h^n$, we see that $A=b^n B$.
\qed

\begin{remark}
In~\cite{MR2679152},
another generalization of B.~Segre's result~\cite[n.~81]{MR0213949} is obtained.
\end{remark}

%
%

\bibliographystyle{plain}

\def\cftil#1{\ifmmode\setbox7\hbox{$\accent"5E#1$}\else
  \setbox7\hbox{\accent"5E#1}\penalty 10000\relax\fi\raise 1\ht7
  \hbox{\lower1.15ex\hbox to 1\wd7{\hss\accent"7E\hss}}\penalty 10000
  \hskip-1\wd7\penalty 10000\box7} \def\cprime{$'$} \def\cprime{$'$}
  \def\cprime{$'$} \def\cprime{$'$}

\end{document}